\documentclass[12pt]{amsart}
\usepackage{amssymb,amsmath,amsthm}
\numberwithin{equation}{section}
\makeatletter
\def\@cite#1#2{{\m@th\upshape\bfseries%
[{#1\if@tempswa{\m@th\upshape\mdseries, #2}\fi}]}}
\makeatother
\newtheorem{thm}{Theorem}[section]

\newtheorem{claim}[thm]{Claim}
\newtheorem{cor}[thm]{Corollary}
\newtheorem{prop}[thm]{Proposition}

\newtheorem{defn}[thm]{Definition}

\newtheorem{ques}[thm]{Question}

\newtheorem{conj}[thm]{Conjecture}

\begin{document}
\title[Norms and spectral radii of composition operators]{Norms and spectral radii of linear fractional composition operators on the ball}
\author{Michael~T.~Jury}
\address{Department of Mathematics,
        University of Florida, 
        Gainesville, Florida 32603}
\email{mjury@math.ufl.edu}
\thanks{Partially supported by NSF grant DMS-0701268}
\date{\today}
\begin{abstract}We give a new proof that every linear fractional map of the unit ball induces a bounded composition operator on the standard scale of Hilbert function spaces on the ball, and obtain norm bounds analogous to the standard one-variable estimates.  We also show that  Cowen's one-variable spectral radius formula extends to these operators.  The key observation underlying these results is that every linear fractional map of the ball belongs to the Schur-Agler class. 
\end{abstract}
\maketitle

\section{Introduction}
\subsection{Background}
Given a set $\Omega$, a collection of functions $\mathcal F :\Omega\to \mathbb C$ and a map $\varphi:\Omega\to \Omega$, one can define a \emph{composition operator}
$$
C_\varphi :f\to  f\circ \varphi
$$
Typically $\Omega$ is a domain in $\mathbb{C}$ or $\mathbb{C}^m$, $\varphi$ is a holomorphic map and $\mathcal F$ is a Banach space of holomorphic functions.  Broadly, one is interested in extracting properties of $C_\varphi$ acting on $\mathcal F$ (boundedness, spectral properties, etc.) from function theoretic or dynamical properties of $\varphi$.     The most studied case is that of $\Omega=\mathbb{D}$, the open unit disk in $\mathbb{C}$, and $\mathcal F$ the Hardy space $H^2$.  In this case it follows from the Littlewood subordination principle that every holomorphic self-map $\varphi$ of $\mathbb{D}$ induces a bounded composition operator on $H^2$.  A theorem of C. Cowen computes the spectral radius of $C_\varphi$.  The purpose of the present paper is to extend Cowen's theorem to a certain class of composition operators acting on the standard scale of holomorphic spaces on the open unit ball $\mathbb{B}^m \subset \mathbb{C}^m$. 

 The primary difficulty in studying composition operators on the ball is that not every holomorphic self map $\varphi$ induces a bounded composition operator on the standard spaces.  Moreover, in many cases even when boundedness can be established, it is difficult to obtain useful norm estimates.  In \cite{jury-pams} we showed that every self-map $\varphi$ of the ball belonging to the \emph{Schur-Agler class} $\mathcal{S}_m$ (defined below) induces a bounded composition operator on the standard scale of spaces, and moreover obeys a norm estimate analogous to the one-variable case.  Since every self-map of the unit disk belongs to the Schur-Agler class, one's intuition is that the maps $\varphi\in\mathcal{S}_m$ should have more behavior in common with self-maps of the disk than do generic self-maps of the ball.

In this paper we show that the linear fractional maps of $\mathbb{B}^m$ introduced by Cowen and MacCluer \cite{MR1768872} belong to the Schur-Agler class and obtain norm bounds.  We then use this result together with an explicit parametrization of the non-elliptic linear fractional maps obtained by Bracci et al.\,\cite{bracci-contreras-diazmadrigal} to obtain a formula for the spectral radius, which extends Cowen's result to linear fractional maps in higher dimensions.  Moreover we conjecture that this formula should hold for all maps in the Schur-Agler class.  

The paper is organized as follows:  we conclude this introductory section by defining the Schur-Agler class $\mathcal S_m$ and describing its relevant properties.  In Section~\ref{S:lft} we prove that every linear fractional map of $\mathbb{B}^m$ belongs to $\mathcal S_m$ and obtain a norm estimate for the induced composition operators; from the norm estimate we deduce a prototype expression for the spectral radius.  In Section~\ref{S:specrad} we prove the spectral radius formula for linear fractional maps and describe some of the geometric difficulties (absent in the one-variable case) encountered in trying to extend the formula to all Schur-Agler mappings.
\subsection{The Schur-Agler class}
Let $\mathbb{B}^m$ denote the open unit ball of $\mathbb{C}^m$.  We will write $\langle \cdot , \cdot \rangle$ for the standard Hermitian inner product on $\mathbb{C}^m$ and $|z|=\sqrt{\langle z,z\rangle}$ for the Euclidean length.  It will often be convenient to write points of $\mathbb C^m$ in the form $z=(z_1, z^\prime)$ with $z_1\in \mathbb C$ and $z^\prime =(z_2, \dots z_m)\in\mathbb C^{m-1}$.
\begin{defn}
The \emph{Schur-Agler class} $\mathcal{S}_m$ is the set of all holomorphic mappings $\varphi:\mathbb{B}^m\to \mathbb{B}^m$ for which the Hermitian kernel
\begin{equation}\label{E:dBR}
k^\varphi(z,w)=\frac{1-\langle \varphi(z), \varphi(w)\rangle}{1-\langle z,w \rangle}
\end{equation}
is positive semidefinite.
\end{defn}
The kernel (\ref{E:dBR}) will be called the \emph{de Branges-Rovnyak kernel} associated to $\varphi$.  When $m=1$ these are the classical de Branges-Rovnyak kernels \cite{MR0215065, MR1289670}.  The functions $\varphi$ for which $k^\varphi$ is positive are precisely those admitting a representation as a transfer function of a multivariate linear system \cite{MR1846055}, but we will not use this representation explicitly. 

It is an elementary but important fact that $\mathcal{S}_m$ is closed under composition:

\begin{thm}
If $\varphi, \psi\in\mathcal{S}_m$ then so is $\varphi\circ\psi$.
\end{thm}
\begin{proof}
The kernel $k^{\varphi\circ\psi}$ may be factored as
$$
k^{\varphi\circ\psi}(z,w)=k^\varphi(\psi(z), \psi(w))\cdot k^\psi(z,w)
$$
which is a pointwise product of positive kernels and hence positive.

\end{proof}
\noindent In particular iterates of Schur-Agler mappings remain in the Schur-Agler class.  It will be proved in the next section that every linear fractional map of $\mathbb{B}^m$ belongs to $\mathcal{S}_m$; in particular every automorphism of the ball belongs to the Schur-Agler class.  
\begin{defn}
Let $m, \beta$ be positive integers.  The space $H^2_{m,\beta}$ is the space of holomorphic functions on the unit ball $\mathbb{B}^m$ with reproducing kernel
$$
k_\beta(z,w)=\frac{1}{(1-\langle z,w\rangle)^\beta}
$$ 
\end{defn}
When $\beta=1$ this is the Drury-Arveson space, which is strictly smaller than the classical Hardy space on the ball but often the more appropriate setting for multivariable operator theory; see e.g. \cite{MR1882259, MR1668582}.  When $\beta=m$ we obtain the classical Hardy space and $\beta=m+1$ gives the Bergman space.  This scale of spaces can be extended to non-integral values of $\beta$ via Calderon interpolation, and all of the results of this paper are valid for this larger scale.  However since the primary values of interest are $\beta=1, m$ and $m+1$, we omit the details.

It was shown in \cite{jury-pams} that every $\varphi\in\mathcal{S}_m$ induces a bounded composition operator on each of the spaces $H^2_{m, \beta}$, satisfying a ``one-variable style'' norm estimate, in particular an estimate which depends only on the value of $\varphi$ at $0$.  In fact when $m=1$ this is precisely the ``classical'' norm estimate for composition operators on the standard scale of Hilbert function spaces.    In higher dimensions, a related upper bound was obtained by Bayart \cite[Theorem 4.1]{MR2296311}, which applies to certain univalent mappings (not necessarily in $\mathcal{S}_m$) but which depends both on $\varphi(0)$ and on global estimates for derivatives of $\varphi$.  
\begin{thm}\label{T:SA_boundedness}
If $\varphi\in\mathcal{S}_m$ then $C_\varphi$ is bounded on $H^2_{m, \beta}$ and 
\begin{equation}\label{E:SA_boundedness}
\left(\frac{1}{1-|\varphi(0)|^2}\right)^{\beta/2}\leq \|C_\varphi\|\leq \left( \frac{1+|\varphi(0)|}{1-|\varphi(0)|}\right)^{\beta/2}
\end{equation}
\end{thm}
\begin{proof}
The upper bound is proved in \cite{jury-pams}; the lower bound is generic for composition operators acting on reproducing kernel Hilbert spaces:  since $k_\beta(\cdot,0)\equiv 1$, 
\begin{align*}
\|C_\varphi^*\| &\geq  \|C_\varphi^* k_\beta(\cdot,0)\|\\ &= \| k_\beta (\cdot,\varphi(0))\| \\&=   \left(\frac{1}{1-|\varphi(0)|^2}\right)^{\beta/2} 
\end{align*}
\end{proof}
We obtain immediately an expression for the spectral radius of $C_\varphi$.  In what follows we let $\varphi_n$ denote the $n^{th}$ iterate of $\varphi$, and observe that $C_\varphi^n =C_{\varphi_n}$.
\begin{cor}\label{C:proto_specrad}
If $\varphi\in\mathcal{S}_m$ then the spectral radius of $C_\varphi$ acting on $H^2_{m,\beta}$ is
\begin{equation}\label{E:proto_specrad}
\lim_{n\to \infty} (1-|\varphi_n(0)|)^{-\beta/2n}
\end{equation}
\end{cor}
\begin{proof}
Since $\mathcal{S}_m$ is closed under composition, we may iterate the norm inequality (\ref{E:SA_boundedness}) to obtain
$$
\|C_\varphi^n\|=\|C_{\varphi_n}\|\sim(1-|\varphi_n(0)|)^{-\beta/2}
$$
Since $r(C_\varphi)=\lim \|C_\varphi^n\|^{1/n}$, the corollary follows.
\end{proof}
The expression (\ref{E:proto_specrad}) should not really be regarded as a formula for the spectral radius, unless some method of evaluating the limit is available.  In one dimension (for maps without interior fixed points), the limit can be evaluated in terms of the angular derivative at the Denjoy-Wolff point.  The evaluation of this limit for linear fractional mappings in higher dimensions is the purpose of the next section; we obtain a result analogous to the one-variable case, where the dilatation coefficient (defined below) plays the role of the angular derivative.

Intuitively, one may expect that Schur-Agler mappings of $\mathbb{B}^m$ may exhibit a stronger affinity with self-maps of $\mathbb D$ than do generic self-maps of $\mathbb{B}^m$.  The reason for this is that \emph{every} self map of $\mathbb{D}$ belongs to $\mathcal S_1$, while for $m>1$ $\mathcal S_m$ is always a proper subset of the self-maps of $\mathbb{B}^m$.  In particular any fact about self-maps of $\mathbb{D}$ which can be proved using only the positivity of the de Branges-Rovnyak kernel ought to have an analogue for the Schur-Agler class; though of course this analogy cannot be taken too literally.  

Finally, a bit of notation:  given two sequences of positive numbers $a_n, b_n$, we write $a_n\sim b_n$ to mean that there exists strictly positive constants $C_1, C_2$ such that 
$$
C_1 \leq \frac{a_n}{b_n}\leq C_2
$$
for all $n$.  

\section{Linear fractional maps}\label{S:lft}
  We now prove that the linear fractional maps of $\mathbb{B}^d$ introduced by Cowen and MacCluer \cite{MR1768872} belong to $\mathcal{S}_m$.  By the theorem and its corollary we obtain a new proof of the boundedness of linear fractional composition operators on the standard spaces, as well as the norm estimate (\ref{E:SA_boundedness}).  

Following Cowen and MacCluer \cite{MR1768872}, a linear fractional map on $\mathbb{B}^m$ is defined to be a function of the form
\begin{equation}
\varphi(z)=\frac{Az+B}{\langle z, C\rangle +D}
\end{equation}
where $A$ is a $m\times m$ matrix, $B,C$ are column vectors in $\mathbb{C}^m$, and $D$ is a complex number.  Here $\langle \cdot, \cdot\rangle$ denotes the standard inner product on $\mathbb{C}^m$.  Clearly, the parameters $A,B, C, D$ are not uniquely determined, since they may all be multiplied by a fixed scalar without changing $\varphi$.  It is shown in \cite{MR1768872} that such map takes $\mathbb{B}^m$ into itself if and only if for some choice of $A,B, C, D$ representing $\varphi$, the $(m+1)\times (m+1)$ matrix
\begin{equation}
T=\begin{pmatrix}
A & B \\ C^* & D 
\end{pmatrix}
\end{equation}
is contractive with respect to the indefinite bilinear form on $\mathbb{C}^{m+1}$ defined by
\begin{equation}
[v,w]=\langle Jv, w\rangle
\end{equation}
where $J$ is the matrix
\begin{equation}
 J= \begin{pmatrix}
I_m & 0 \\ 0 & -1 
\end{pmatrix}
\end{equation}
That is, $T$ must satisfy
$$
[Tv, Tv]\leq [v,v]
$$
for all $v\in\mathbb{C}^{m+1}$.  This contractivity condition is satisfied if and only if the matrix $J-T^* JT$ is positive semidefinite.  We will make use of the condition in this latter form.  

It is then proved in \cite{MR1768872} that every such map induces a bounded composition operator on the standard scale of spaces (at least when $\beta \geq m$), though this proof is indirect and in particular does not provide an estimate for the norm of $C_\varphi$.   We will prove that $C_\varphi$ is bounded by appeal to Theorem~\ref{T:SA_boundedness}, and prove that the de Branges-Rovnyak kernel $k^\varphi$ is positive by exhibiting an explicit factorization, which we obtain from a factorization of the (assumed positive) matrix $J-T^*JT$.  
We can now state the factorization result:
\begin{thm}\label{T:LFT_factorization}
Every linear fractional map $\varphi:\mathbb{B}^m\to \mathbb{B}^m$ belongs to the Schur-Agler class $\mathcal S_m$.
\end{thm}
\begin{proof}
Let $T$ be a $(m+1)\times(m+1)$ matrix which is contractive with respect to $[\cdot, \cdot]$ and has the form
\begin{equation}
T=\begin{pmatrix}
A & B  \\ C^* & D 
\end{pmatrix}
\end{equation}
and let $\varphi$ denote the associated linear fractional transformation. (By the remarks preceding the proof, every linear fractional self-map of $\mathbb B^m$ arises in this way.)
Factor $J-T^*JT$ as
\begin{equation}
J-T^*JT =X^*X
\end{equation}
with
\begin{equation}
X=\begin{pmatrix}
X_{11} & X_{12} \\ X_{21}^* & X_{22}
\end{pmatrix}
\end{equation}
Now define a function $L:\mathbb{B}^m\to \mathbb{C}^{m+1}$ by
\begin{equation}
L(z) = X\begin{pmatrix} z\\ 1\end{pmatrix}=\begin{pmatrix} X_{11} z +X_{12} \\ \langle z, X_{21}\rangle +X_{22}\end{pmatrix}
\end{equation}
We now claim that the de Branges-Rovnyak kernel can be factored as
\begin{equation}\label{E:factorization}
k^\varphi(z,w)= \frac{1}{\langle z, C\rangle+D} \left(1+\frac{L(z)L(w)^*}{1-\langle z,w\rangle}\right) \frac{1}{\overline{\langle w, C\rangle+D}}
\end{equation}
from which it is apparent that $k^\varphi$ is positive.
To verify (\ref{E:factorization}), we first write out $k^\varphi(z,w)$ as
\begin{align}
\begin{split}
k^\varphi(z,w) &=\frac{1}{\langle z, C\rangle+D} \frac{1}{\overline{\langle w, C\rangle+D}} \\ &\phantom{==}\times \frac{(\langle z, C\rangle+D)\overline{(\langle w, C\rangle+D)} -\langle Az +B, Aw+B\rangle }{1-\langle z,w\rangle}
\end{split}
\end{align}
Working with the factor on the second line, we verify that its numerator is equal to $1-\langle z,w\rangle +L(z)L(w)^*$, which proves (\ref{E:factorization}):
\begin{align*}
1-\langle z,w\rangle +L(z)&L(w)^* = 1-\langle z,w\rangle +\langle X^*X \begin{pmatrix} z \\ 1\end{pmatrix},\begin{pmatrix} w \\ 1\end{pmatrix} \rangle\\
 &= 1-\langle z,w\rangle +\langle J-T^*JT \begin{pmatrix} z \\ 1\end{pmatrix},\begin{pmatrix} w \\ 1\end{pmatrix}\rangle \\ 
 &= -\langle JT \begin{pmatrix} z \\ 1\end{pmatrix},T\begin{pmatrix} w \\ 1\end{pmatrix}\rangle \\
 &= (\langle z, C\rangle+D)\overline{(\langle w, C\rangle+D)} -\langle Az +B, Aw+B\rangle 
\end{align*}

\end{proof}

\section{Spectral radii}\label{S:specrad}

We begin with some basic definitions and results about the iteration of self-maps of the ball, and then describe some known results.  Suppose that $\varphi:\mathbb{B}^m\to \mathbb{B}^m$ is a holomorphic mapping which does not fix any point of $\mathbb{B}^m$.  MacCluer \cite{MR694933} showed that an analogue of the Denjoy-Wolff theorem holds:  there exists a unique point $\zeta\in\partial\mathbb{B}^m$ such that the iterates of $\varphi$ converge uniformly to $\zeta$ on compact subsets of $\mathbb{B}^m$.  This point will be called the \emph{Denjoy-Wolff} point of $\varphi$.  Moreover, it follows from  \cite[Theorem 1.3]{MR694933} that 
$$
0< \liminf_{z\to \zeta} \frac{1-|\varphi(z)|^2}{1-|z|^2} =\alpha \leq 1
$$
and hence by the Julia-Caratheodory theorem on the ball \cite[Theorem 8.5.6]{MR601594} the complex directional derivative $D_\zeta \varphi$ has a radial limit\footnote{In fact this limit exists in the wider sense of \emph{restricted $K$-limit} (or \emph{hypoadmissible limit}) but we will not require this notion at the moment.} $\alpha$ at $\zeta$; this number is called the \emph{dilatation coefficient} of $\varphi$.  (When $m=1$, $\alpha$ is the angular derivative of $\varphi$ at $\zeta$.)  
The following is then a special case of Julia's theorem on the ball (\cite[Theorem 1.3]{MR694933} and \cite[Theorem 8.5.3]{MR601594}):
\begin{thm}
Let $\varphi:\mathbb{B}^m\to \mathbb{B}^m$ with Denjoy-Wolff point $\zeta\in\partial\mathbb B^m$ and dilatation coefficient $\alpha$.  Then for all $z\in \mathbb{B}^m$, 
\begin{equation}\label{E:julia}
\frac{|1-\langle \varphi(z), \zeta\rangle|^2}{1-|\varphi(z)|^2} \leq \alpha\frac{|1-\langle z, \zeta\rangle|^2}{1-|z|^2}
\end{equation}

\end{thm}

We now divide the self-maps of $\mathbb{B}^m$ into three classes:
\begin{defn} A holomorphic self-map $\varphi$ of $\mathbb{B}^m$ will be called
\begin{itemize}
\item \emph{elliptic} if $\varphi$ fixes a point of $\mathbb{B}^m$,
\item \emph{parabolic} if $\varphi$ has no fixed point and dilatation coefficient $1$, and
\item \emph{hyperbolic} if $\varphi$ has no fixed point and dilatation coefficient $\alpha<1$.
\end{itemize}
\end{defn}

In one dimension, Cowen \cite{MR695941} obtained the following formula for the spectral radius of composition operators on $H^2(\mathbb{D})$:
\begin{thm}
Let $\varphi:\mathbb{D}\to \mathbb{D}$.  If $\varphi$ is elliptic then the spectral radius of $C_\varphi$ is $1$; if $\varphi$ is non-elliptic then the spectral radius is $\alpha^{-1/2}$.
\end{thm}
For linear fractional maps in higher dimensions, MacCluer \cite{MR760878} obtained the full spectrum for automorphic symbols $\varphi$ acting on the Hardy space (our case $\beta=m$); it follows from these results that the spectral radius is $1$ for elliptic automorphisms and $\alpha^{-m/2}$ otherwise.   More recently Bayart \cite{MR2296311} obtained the full spectrum for certain parabolic maps conjugate to generalized Heisenberg translations of the Siegel half-space; for these parabolic maps the spectral radius is $1$. 

The spectral radius formulae we obtain will be valid for all elliptic and parabolic maps in the Schur-Agler class; it is only in the hyperbolic case that we restrict to linear fractional maps.  Indeed in the elliptic and parabolic cases the proof we now give is identical to Cowen's in dimension $1$.

\begin{thm}\label{T:spectral_radius} Let $\varphi\in\mathcal{S}_m$.  If $\varphi$ is elliptic or parabolic, then the spectral radius of $C_\varphi$ on $H^2_{d,\beta}$ is $1$.
\end{thm}
\begin{proof}  If $\varphi$ is elliptic, then $C_\varphi$ is similar (via conjugation by an automorphism) to a composition operator $C_\psi$ with $\psi\in\mathcal{S}_m$ and $\psi(0)=0$.  Since $\mathcal{S}_m$ is automorphism invariant, $\psi\in\mathcal{S}_m$ and hence $\|C_{\psi_n}\|=1$ for all $n$ by Theorem~\ref{T:SA_boundedness}, and thus $r(C_\varphi)=r(C_\psi)=1$.  

Now assume $\varphi$ is parabolic with Denjoy-Wolff point $\zeta\in\partial\mathbb{B}^m$.  If $z_n$ is a sequence in $\mathbb{B}^m$ such that $z_n\to \zeta$,  $\varphi(z_n)\to \zeta$, and the limit 
$$
M=\lim_{n\to \infty} \left( \frac{1-|\varphi(z_n)|}{1-|z_n|}\right)
$$ 
exists, then $M\geq 1$.  It follows that 
$$
\liminf_{n\to\infty} \left(\frac{1-|\varphi_n(0)|}{1-|\varphi_{n-1}(0)|}\right)\geq 1
$$
Therefore
\begin{align*}
\lim_{n\to\infty} (1-|\varphi_n(0)|)^{-1/2n} &= \lim_{n\to\infty} \left( \prod_{k=0}^{n-1} \frac{1-|\varphi_k(0)|}{1-|\varphi_{k-1}(0)|}  \right)^{1/2n} \\ & \leq \limsup_{n\to\infty} \left( \frac{1-|\varphi_{n-1}(0)|}{1-|\varphi_n(0)|}\right)^{1/2} \\ &\leq 1
\end{align*}
Thus $r(C_\varphi)\leq 1$ by Corollary~\ref{C:proto_specrad}, and since $1$ is an eigenvalue  $r(C_\varphi)=1$.
\end{proof}
The evaluation of the limit (\ref{E:proto_specrad}) in the hyperbolic case requires a more detailed analysis of the orbit $\{\varphi_n(0)\}$, which can be carried out explicitly in the case of linear fractional maps.  The proof exploits a parametrization of non-elliptic linear fractional maps (conjugated to the Siegel half-space) obtained by Bracci, Contreras and Diaz-Madrigal \cite[Lemma 4.1 and Proposition 4.2]{bracci-contreras-diazmadrigal}.

\begin{thm}\label{T:hyperbolic_specrad}
Let $\varphi$ be a hyperbolic linear fractional map of $\mathbb{B}^m$ with dilatation coefficient $\alpha<1$.  Then
\begin{equation}\label{E:hyperbolic_specrad}
\lim_{n\to \infty} (1-|\varphi_n(0)|^2)^{1/n} =\alpha
\end{equation}
\end{thm}

\begin{proof}
Conjugating $\varphi$ by a rotation of $\mathbb C^m$, we may assume the Denjoy-Wolff point is $e_1=(1, 0,\dots 0)$; clearly (\ref{E:hyperbolic_specrad}) is unchanged.

It will be convenient to move the problem to the Siegel right half-space
$$
\mathbb{H}^m =\{(w_1,w^\prime)\in\mathbb{C}\times\mathbb{C}^{m-1}:\text{Re}\ w_1>\|w^\prime\|^2\}
$$
which is biholomorphically equivalent to $\mathbb{B}^m$ via the generalized Cayley transform
$$
\psi(z_1, z^\prime)=\left( \frac{1+z_1}{1-z_1}, \frac{z^\prime}{1-z_1}\right)
$$
and its inverse
$$
\psi^{-1}(w_1, w^\prime)=\left(\frac{w_1-1}{w_1+1}, \frac{2w^\prime}{w_1+1}\right)
$$
This correspondence extends continuously to identify $\partial\mathbb{B}^m$ with the one-point compactification of $\partial\mathbb{H}^m$, with $e_1$ taken to the point at infinity.

In particular one may calculate that for any $z=(z_1, z^\prime)\in\mathbb{B}^m$, if $w=\psi(z)$ then
$$
1-|z|^2= \frac{4}{|w_1+1|^2}(\text{Re}\ w_1 -\|w^\prime\|^2)
$$
By \cite[Lemma 4.1]{bracci-contreras-diazmadrigal} a map $\varphi$ satisfying our hypotheses is conjugate to a map $\tilde{\varphi}:\mathbb{H}^n \to \mathbb{H}^n$ of the form
\begin{equation}\label{E:BCD_repn}
\tilde\varphi(w_1, w^\prime)=\frac1\alpha (w_1+c+\langle w^\prime,b\rangle , Aw^\prime+d\rangle
\end{equation}
for suitable scalar $c\in\mathbb{C}$, vectors $b,d\in\mathbb{C}^{m-1}$ and $(m-1)\times (m-1)$ matrix $A$.  Of course these parameters satisfy a number of relations, determined by the condition that $\tilde\varphi$ maps $\mathbb H^m$ into itself; the only one we will require explicitly is the fact that $\|A\|\leq \alpha^{1/2} <1$ \cite[Lemma 4.1(i)]{bracci-contreras-diazmadrigal}.  
Let us now write
$$
\tilde\varphi_n(1,0)=(u_n, v_n)
$$
with $u_n\in\mathbb{C}, v_n\in\mathbb{C}^{m-1}$.  Our goal is now to show that
\begin{equation}\label{E:halfspace_specrad}
\lim_{n\to \infty} \left( \frac{4}{|u_n+1|^2}(\text{Re}\ u_n -\|v_n\|^2)\right)^{1/n}=\alpha
\end{equation}
Since $\tilde{\varphi}$ has Denjoy-Wolff point $\infty$, it follows in particular that $|u_n|\to \infty$ and hence $|u_n|\sim |u_n +1|$.  Thus, to establish (\ref{E:halfspace_specrad}) it suffices to show
\begin{equation}
|u_n|\sim \frac{1}{\alpha^n}
\end{equation}
and
\begin{equation}
(\text{Re}\ u_n -\|v_n\|^2)\sim \frac{1}{\alpha^n}
\end{equation}

To do this, we will obtain fairly explicit expressions for $u_n$ and $v_n$; we begin by introducing some notation.  
For each integer $n\geq 0$ define
$$
\beta_n=\sum_{k=0}^n \alpha^k
$$
and polynomials
$$
p_n(z)=\sum_{k=0}^n \beta_{n-k} z^k, \qquad q_n(z)=\sum_{k=0}^n \alpha^{n-k} z^k 
$$
It is straightforward to verify the following recurrence relations:
\begin{gather}
\beta_{n+1}=\alpha \beta_n +1 \\
p_{n+1}(z)=\alpha p_n(z) +\sum_{k=0}^{n+1} z^k \\
q_{n+1}(z)=zq_n(z)+\alpha^n 
\end{gather}
Using these one may also deduce
\begin{equation}
q_n(z)+p_{n-1}(z)=p_n(z)
\end{equation}
With these identities established one can verify by induction that
$$
\tilde{\varphi}(1,0)=\frac1\alpha(1+c, d)
$$
and for all $n\geq 2$
$$
\tilde{\varphi}_n(1,0)=\frac{1}{\alpha^n}\left(1+\beta_{n-1} c +\langle p_{n-2}(A)d, b\rangle, q_{n-1}(A)d\right)
$$
So in particular
$$
u_n = \frac{1}{\alpha^n}(1+\beta_{n-1} c +\langle p_{n-2}(A)d, b\rangle).
$$
Now define
$$
x_n:=\alpha^n u_n = (1+\beta_{n-1} c +\langle p_{n-2}(A)d, b\rangle).
$$
We observe that the real part of $x_n$ must always be strictly positive, and we will show that $x_n\to x$ with $\text{Re}\ x\geq 1$.  This establishes the claimed asymptotic behavior of $|u_n|$.  

The convergence of $x_n$ depends upon the convergence of the polynomials $p_n$; in particular the following fact:
\begin{claim}
The sequence of polynomials $p_n$ converges to 
$$
\frac{1}{1-\alpha}\frac{1}{1-z}
$$ 
uniformly in the disk $|z|\leq \sqrt{\alpha}$. 
\end{claim}
\emph{Proof of claim:}  Let $\|\cdot \|_\infty$ denote the supremum norm over the closed disk of radius $\sqrt{\alpha}$.  Then for every $n$
\begin{align}
\left\|(1-\alpha)p_n -\sum_{k=0}^{n+1} z^k\right\|_\infty &=\left\| \sum_{k=0}^{n+1} \alpha^{n-k+1}z^k\right\|_\infty \\
 &\leq \alpha^{n+1} \left\|\sum_{k=0}^{n+1} \alpha^{-k}z^k\right\|_\infty \\
&\leq \alpha^{n+1}\frac{\alpha^{-(n+2)/{2}}-1}{\alpha^{-1/2}-1}
\end{align}
which tends to $0$ as $n\to \infty$.  Since $\sum_{k=0}^n z^k\to (1-z)^{-1} $ uniformly in this disk, the claim is proved.

Using now the crucial fact that $\|A\|\leq \sqrt{\alpha}$, we conclude that $x_n$ converges to 
\begin{equation}\label{E:xdef}
x=1+\frac{1}{1-\alpha}(c+\langle (I-A)^{-1}d, b\rangle )
\end{equation}
 Now define 
$$
u=(I-A)^{-1}d
$$
and observe that $Au+d=u$.  Since $\tilde{\varphi}$ maps the closure of $\mathbb{H}^n$ into itself, it follows that $\tilde\varphi(\|w^\prime\|^2, w^\prime)\in\overline{\mathbb H^m}$ for all $w^\prime\in \mathbb C^{m-1}$; that is, 
$$
\alpha \|w^\prime\|^2 +\alpha \text{Re}\ \langle w^\prime,b\rangle +\alpha \text{Re}\ c \geq \|Aw^\prime+d\|^2.
$$
Applying this with $w^\prime=u$ gives
$$
\text{Re}\ \langle u,b\rangle +\text{Re}\ c \geq \frac{1-\alpha}{\alpha} \|u\|^2 \geq 0
$$
and hence $\text{Re}\ x\geq 1$.

We now consider $\alpha^n (\text{Re}\ u_n-\|v_n\|^2)$.  Since $\text{Re}\ u_n-\|v_n\|^2\geq 0$, the upper bound follows immediately from the upper bound for $\alpha^n|u_n|$.  To prove boundedness from below, we return momentarily to the ball.  By induction on Julia's theorem (\ref{E:julia}),
$$
\frac{|1-\langle \varphi_n(0), e_1\rangle|^2}{1-|\varphi_n(0)|^2}\leq \alpha^n
$$
for all $n$.  Transferring this inequality to $\mathbb{H}^m$ we obtain
$$
\alpha^n (\text{Re}\ u_n-\|v_n\|^2)\geq 1
$$
for all $n$.  

\end{proof}

\begin{cor}
If $\varphi$ is a hyperbolic linear fractional map of $\mathbb{B}^m$ with dilatation coefficient $\alpha$, then the spectral radius of $C_\varphi$ acting on $H^2_{m, \beta}$ is $\alpha^{-\beta/2}$.
\end{cor}
\begin{proof}
Combine Theorem~\ref{T:hyperbolic_specrad} and Corollary~\ref{C:proto_specrad}.
\end{proof}

To summarize, combining the two spectral radius results we have extended Cowen's spectral radius formula to linear fractional maps in higher dimensions:
\begin{thm}\label{T:lft_spectral_radius}
Let $\varphi$ be a linear fractional self-map of $\mathbb{B}^m$.  The spectral radius of $C_\varphi$ acting on $H^2_{m,\beta}$ is $1$ if $\varphi$ is elliptic; if $\varphi$ is non-elliptic with dilatation coefficient $\alpha$ the spectral radius is $\alpha^{-\beta/2}$.
\end{thm}

\begin{conj}
The spectral radius formulae of Theorem~\ref{T:lft_spectral_radius} are valid for all $\varphi\in\mathcal{S}_m$.   
\end{conj}

\noindent By Theorem~\ref{T:spectral_radius} the conjecture is true for elliptic and parabolic maps.

In the hyperbolic case, one may try to prove the conjecture by a method analogous to Cowen's proof in the disk \cite{MR695941}, namely, by proving that the iterates $\varphi_n(0)$ converge to the Denjoy-Wolff point sufficiently well so that
$$
\lim_{n\to \infty}\frac{1-|\varphi_n(0)|^2}{1-|\varphi_{n-1}(0)|^2}=\alpha
$$
In one variable, this is accomplished by showing that when $\alpha <1$, the iterates $\varphi_n(0)$ converge nontantgentially to the Denjoy-Wolff point; the above limit then follows from the Julia-Caratheodory theorem.  In the ball, one needs \emph{restricted} convergence in order to invoke the corresponding version of Julia-Caratheodory:  to define this, fix a point $\zeta\in\partial\mathbb{B}^n$ and consider a curve $\Gamma:[0, 1)\to \mathbb{B}^n$ such that $\Gamma(t)\to \zeta$ as $t\to 1$.  Let $\gamma(t)=\langle \Gamma(t), \zeta\rangle\zeta$ be the projection of $\Gamma$ onto the complex line through $\zeta$.  The curve $\Gamma$ is called \emph{special} if 
\begin{equation}
\lim_{t\to 1}\frac{|\Gamma -\gamma|^2}{1-|\gamma|^2}=0
\end{equation}
and \emph{restricted} if it is special and in addition
\begin{equation}
\frac{|\zeta-\gamma|}{1-|\gamma|^2}\leq A
\end{equation}
for some constant $A>0$.  We say that a function $f:\mathbb{B}^n\to \mathbb{C}$ has \emph{restricted $K$-limit} $L$ at $\zeta$ if $\lim_{z\to \zeta}f(z)=L$ along every restricted curve.  Now, if $\varphi$ is a non-elliptic self-map of $\mathbb{B}^m$ with Denjoy-Wolff point $\zeta$ and dilatation coefficient $\alpha$, it follows from the Julia-Caratheodory theorem that the function
$$
\frac{1-|\varphi(z)|^2}{1-|z|^2}
$$
has restricted $K$-limit $\alpha$ at $\zeta$.  Thus the conjecture is true for any hyperbolic $\varphi$ for which $\varphi_n(0)\to \zeta$ restrictedly.  However the following shows that in general we need not have restricted convergence, even for linear fractional maps.

\begin{prop}
Let $\varphi$ be a hyperbolic linear fractional map with with Denjoy-Wolff point $e_1$ and dilatation coefficient $\alpha$, and let $\tilde\varphi$ be the conjugate mapping of $\mathbb{H}^m$ given by (\ref{E:BCD_repn}).  If $\varphi_n(0)\to e_1$ restrictedly, then 
\begin{equation}\label{E:restricted_little_o}
\|q_{n-1}(A)d\|^2=o(\alpha^n).
\end{equation}

\end{prop}
\begin{proof}
If $\varphi_n(0)\to e_1$ restrictedly then 
$$
\lim_{n\to \infty} \frac{|\varphi_n(0)-\langle \varphi_n(0), e_1\rangle |^2}{1-|\langle \varphi_n(0), e_1\rangle |^2} =0
$$
Under the Cayley transform, this is equivalent to 
$$
\lim_{n\to \infty} \frac{\|v_n^2\|}{\text{Re}\ u_n} =0
$$
which is in turn the same as 
$$
\lim_{n\to \infty} \frac{1}{\alpha^n}\frac{\|q_{n-1}(A)d\|^2}{\text{Re}\ x_n}=0
$$
Since $\text{Re}\ x_n\sim 1$, this proves the theorem.
\end{proof}
Using the parametrization (\ref{E:BCD_repn}) it is straightforward to construct hyperbolic linear fractional maps for which the condition (\ref{E:restricted_little_o}) does not hold.\footnote{The corresponding ``big O'' condition is always satisfied.}  To do this, fix $0<\alpha<1$ and let $A$ be the diagonal matrix with each diagonal entry equal to $\sqrt{\alpha}$.  Let $d$ be any unit vector in $\mathbb{C}^{m-1}$ and define $b=2\alpha^{-1/2} d$, $c=\alpha^{-1}$.  Then $\tilde\varphi$ defined by (\ref{E:BCD_repn}) is a conjugate to a hyperbolic linear fractional map for which (\ref{E:restricted_little_o}) is violated:  we calculate
$$
\alpha^{-n}\|q_{n-1}(A)d\|^2 = \alpha^n \left( \sum_{k=0}^n \alpha^{-k/2}\right)^2 =\left( \frac{1-\alpha^{(n+1)/2}}{1-\alpha^{1/2}}\right)^2
$$
which is greater than $1$ for all $n$.

Even though the orbit $\varphi_n(0)$ need not approach the Denjoy-Wolff point restrictedly, it can be shown (at least when $m=2$) that when $\varphi$ is a linear fractional map, the limit
$$
\lim_{n\to \infty}\frac{1-|\varphi_n(0)|^2}{1-|\varphi_{n-1}(0)|^2}
$$
exists and equals $\alpha$.  We do not know if this is true of general Schur-Agler mappings.

\begin{ques}
If $\varphi\in\mathcal{S}_m$ is hyperbolic with dilatation coefficient $\alpha$, is it true that
$$
\lim_{n\to \infty}\frac{1-|\varphi_n(0)|^2}{1-|\varphi_{n-1}(0)|^2}
$$
exists and equals $\alpha$?
\end{ques}
An affirmative answer to this question would prove the conjecture.  If on the other hand the limit exists for some $\varphi$ but has a value different from $\alpha$ (necessarily larger) then the conjecture would be false.  One may try to answer the question by looking for a stronger form of the Julia-Caratheodory theorem in the ball (valid for Schur-Agler mappings).  Some results in this direction are obtained in \cite{jury-preprint}, but so far these results are not sufficient to answer the question.
\bibliographystyle{plain} 
\bibliography{linfrac} 
\end{document}